\documentclass[12pt]{article}

\usepackage[english]{babel}

\usepackage[T1]{fontenc}

\usepackage[a4paper,top=1in,bottom=1in,left=1in,right=1in,marginparwidth=1.75cm]{geometry}

\usepackage{amsmath}

\usepackage[colorinlistoftodos]{todonotes}

\usepackage{amsfonts}
\usepackage{xstring}
\usepackage{csquotes}
\usepackage{biblatex}
\DeclareLanguageMapping{english}{english-apa}
\usepackage{titleps}
\usepackage{csquotes}
\usepackage{amssymb}
\usepackage{amsthm}
\usepackage{tikz}
\usepackage{setspace}

\let\oldaddcontentsline\addcontentsline
\newcommand{\stoptocentries}{\renewcommand{\addcontentsline}[3]{}}
\newcommand{\starttocentries}{\let\addcontentsline\oldaddcontentsline}

\newlength{\normalparindent}
\raggedright
\usetikzlibrary{arrows}
\title{Mathieu-Zhao Subspaces of Vertex Algebras}
\author{Matthew Speck}
\date{}

\newcommand\leftidx[3]{%
  {\vphantom{#2}}#1#2#3%
}
\theoremstyle{definition}
\newtheorem{theorem}{Theorem}
\newtheorem{definition}[theorem]{Definition}
\newtheorem{corollary}[theorem]{Corollary}
\newtheorem{lemma}[theorem]{Lemma}
\newtheorem{conjecture}[theorem]{Conjecture}
\newtheorem{example}[theorem]{Example}

\newtheorem{proposition}[theorem]{Proposition}
 
\theoremstyle{remark}

\begin{document}
\begin{center}

MATHIEU-ZHAO SUBSPACES OF VERTEX ALGEBRAS\\
\vspace{1cm}
MATTHEW SPECK\\

\end{center}
\noindent 31 Pages\\

A Mathieu-Zhao subspace is a generalization of an ideal of an associative algebra $\mathcal A$ over a unital ring $R$ first formalized in 2010.  A vertex algebra is an algebraic structure first developed in conjunction with string theory in the 1960s and later axiomatized by mathematicians in the 1980s.  We formally introduce the definition of a Mathieu-Zhao subspace $M$ of a vertex algebra $V$.  Building on natural connections to associative algebras, we classify an infinite set of non-trivial, non-ideal Mathieu-Zhao subspaces for simple and general vertex algebras by group action eigenspace decomposition.  Finally, we state the locally nilpotent $\varepsilon$-derivation (LNED) conjecture for vertex algebras.\\
\vspace*{1cm}
\noindent
KEYWORDS: Vertex algebra, Ideal, Mathieu-Zhao subspace, Eigenspace decomposition \newpage
\begin{titlepage}
    \begin{center}
       
        MATHIEU-ZHAO SUBSPACES OF VERTEX ALGEBRAS\\
		\vspace{1cm}
		MATTHEW SPECK\\
        
        \vfill
        \singlespacing
        A Thesis Submitted in Partial\\
        Fulfillment of the Requirements\\
        for the Degree of\\ \doublespacing
        MASTER OF SCIENCE\\
        
        Department of Mathematics\\
        ILLINOIS STATE UNIVERSITY\\
        2018
        
    \end{center}
\end{titlepage}

\begin{center}
Copyright 2018 Matthew Speck
\end{center}\newpage
\begin{center}
       
       MATHIEU-ZHAO SUBSPACES OF VERTEX ALGEBRAS\\
		\vspace{1cm}
		MATTHEW SPECK\\
\end{center}
\vfill
\begin{flushright}
        $\begin{array}{l}
        \text{COMMITTEE MEMBERS:} \\
        \text{Gaywalee Yamskulna, Chair}\\
        \text{Fusun Akman} \\
        \text{Wenhua Zhao}
        \end{array}$
\end{flushright}
\newpage
\pagenumbering{roman}
\addcontentsline{toc}{section}{Acknowledgments}
\begin{center}
ACKNOWLEDGMENTS
\end{center}
\begin{flushleft}
My deepest and sincerest thanks to the following who have motivated me and maintained faith in me over the course of my time at Illinois State University.\\
\vspace{.5cm}
To my supervisor, Gaywalee Yamskulna, whose patience, understanding, and excitement for all of her students constantly inspires me.\\
\vspace{.5cm}
To my other committee members, Fusun Akman and Wenhua Zhao, whose lessons on professionalism, rigor, and how to drive to Chicago will be with me throughout my career.\\
\vspace{.5cm}
To Alberto Delgado for giving me permission to love math again.\\
\vspace{.5cm}
To the rest of the amazing faculty and staff in the Department of Mathematics who have taught me the true value of a caring and invested education.\\
\vspace{.5cm}
And to my friends and family, without whom I would certainly never have completed this thesis.
\end{flushleft}
\begin{flushright}
M. S.
\end{flushright}
\newpage
\addcontentsline{toc}{section}{Contents}
\begin{center}
CONTENTS
\end{center}
\hfill Page\\
ACKNOWLEDGMENTS\hfill i\\
CONTENTS\hfill ii\\
CHAPTER I: MATHIEU-ZHAO SUBSPACES\hfill 1\\

\hspace{.5cm}History\hfill 1\\

\hspace{.5cm}Generalizing Ideals\hfill 4\\
CHAPTER II: VERTEX ALGEBRAS AND THEIR MATHIEU-ZHAO SUBSPACES\hfill 11\\
\hspace{.5cm}Vertex Algebras\hfill 11\\
\hspace{.5cm}Main Results\hfill 18\\
\hspace{.5cm}Next Steps\hfill 26\\
REFERENCES\hfill 30

\stepcounter{section}

\newpage
\pagenumbering{arabic}
\setcounter{section}{1}
\stoptocentries
\begin{center}
CHAPTER I: MATHIEU-ZHAO SUBSPACES
\end{center}

\begin{center}
\textbf{History}
\end{center}

First axiomatized by Richard Borcherds in the 1980s in his (ultimately successful and Field's Medal-worthy) attempt at proving the famous Moonshine Conjecture (\cite{moonshine}), vertex algebras saw their first iterations in the early days of string theory.  The quantum physicists working on the project needed a way to describe particle and string interactions, and their new theory focused on the point of interaction.  In the spirit of Feynman diagrams, this moment was illustrated as a vertex, and as we will see, vertex algebras are a study of the actual operations which occur at a given vertex.  As such, a vertex algebra can be understood as a vector space of functions on an underlying vector space, in much the same way as the polynomial ring $\mathbb R[x]$ is a vector space of functions (polynomials) on the underlying vector space $\mathbb R$.  Vertex algebras only appear to be more complicated because of the lack of familiarity we (humble, yet eager mathematicians) have with them and because they lack certain properties that we typically take for granted---most notably, associativity.

The mathematical interest in vertex algebras grew out of questions in group theory.  Particularly, mathematicians needed to classify the last of the sporadic finite simple groups, and the biggest of these was the monster group.  The so-called monstrous moonshine arose as a vertex algebra whose symmetries were defined by the monster group.  The key players in the development of this theory---Borcherds, Conway, Frenkel, Lepowski, McKay, and Meurman---worked together and separately, but each established a far-reaching theory that impacted the entire world of group algebra.  Here was a new object, and we needed to know everything.

In the ever-expanding study of vertex algebras, the most natural question we must ask is: What can we learn of vertex algebras from more established theories?  For example, vertex algebras have a natural connection with Lie algebras which serve, in a sense, as their foundation.  As such, much of the established theory (e.g. \cite{dongintro,zhu1}) capitalizes on this connection quite effectively.  However, what is missed in this angle of development is the sense that Borcherds gave us in (\cite{borcherdstalk}): that vertex algebras are a kind of commutative ring.  With that in mind, what can we learn of vertex algebras by studying commutative rings?

This paper on vertex algebras attempts to define and establish the theory for one concept which has only recently been formalized in the language of (commutative) rings.  Mathieu-Zhao subspaces were first formalized in (\cite{zhaogen}) in the interest of solving the famous Jacobian Conjecture.

\theoremstyle{definition}

\begin{conjecture}[Jacobian Conjecture]
Let $k$ be a field, let $$F:k^N\to k^N,$$
$$F(c_1,\dots, c_N) = (f_1(c_1,\dots, c_N),\dots, f_N(c_1,\dots,c_N))$$ be a function (with polynomial components $f_1,\dots,f_N$), and let $J_F$ denote the Jacobian of $F$.  If $J_F$ is a nonzero constant and if $k$ has characteristic 0, then $F$ has an inverse function $F^{-1} = G:k^N\to k^N$ whose components are also polynomials.
\end{conjecture}

Zhao showed that the Jacobian Conjecture is a special case of the so-called Image Conjecture, which can be stated in terms of Mathieu-Zhao subspaces (\cite{zhaogen}).  For a more detailed background on the Image Conjecture, see (\cite{essen}).

\begin{conjecture}[Image Conjecture]
Define $\mathbb{D}[z]$ to be the set of differential operators $\Phi$ of $\mathbb{C}[z]$ such that $\Phi=h(z)+\sum_{i=1}^n c_i \partial/\partial z_i$ for some $h(z)\in\mathbb{C}[z]$ and $c_i\in\mathbb{C}$.  For any subset $\mathcal{C}=\{\Phi_i|i\in I\}$ of differential operators indexed by a set $I\subset \mathbb Z$, set $\text{Im }\mathcal{C} =\sum_{i\in I}(\Phi_i\mathbb{C}[z])$.  We say that $\mathcal{C}$ is \textit{commuting} if for any $i,j\in I$, $\Phi_i$ and $\Phi_j$ commute.  We conjecture that for any commuting subset $\mathcal C\subset \mathbb D [z]$, $\text{Im }\mathcal C$ is a Mathieu-Zhao subspace of $\mathbb C [z]$.
\end{conjecture}

We will not prove the Jacobian Conjecture or this generalization.  Our target will be to establish a theory of Mathieu-Zhao subspaces on vertex algebras with the specific goal of developing a basic understanding of the LNED conjecture (\cite{zhaoopen}) on vertex algebras.  These conjectures are already being studied for the (associative) polynomial ring $K[x]$ over a field $K$.  Throughout the course of this paper, we will see intimate connections develop between this polynomial ring and a given vertex algebra.  Higher-dimensional rings will correspond to higher-dimensional vertex algebras, and conjectures in rings will have natural analogues in corresponding VA's.  Suffice it to say, we will not prove any of these major conjectures, but instead, our goal will be to establish relevant analogies between these rings and their corresponding VA's, state these conjectures in this new context, and hope that the insights provided therein will lead to new lenses through which to attack these extremely stubborn algebraic quandaries.

\begin{definition}[\cite{zhaoopen}]
Let $\mathcal A$ be an algebra over a unital, associative ring $R$.  An \textit{$R$-derivation} $D$ is an $R$-linear map
$$D:\mathcal A\to\mathcal A$$
such that
$$D(ab)=D(a)b+aD(b),$$
and an \textit{$R$-$\varepsilon$-derivation} $\partial$ is an $R$-linear map
$$\delta:\mathcal A\to \mathcal A$$
such that
$$\delta(ab)=b\delta(a)+a\delta(b)-\delta(a)\delta(b).$$
A derivation $\partial$ (of either type) is called \textit{locally finite} if for each $a\in\mathcal A$, the $R$-submodule spanned by $\partial^i(a)$ (with $i\geq 0$) over $R$ is finitely generated, and $\partial$ is called \textit{locally nilpotent} if for each $a\in\mathcal A$, there exists $m\geq 1$ such that $\partial^m(a)=0$.
\end{definition}

Understanding these derivations, we can now state the following conjecture, called the locally nilpotent $\varepsilon$-derivation conjecture (LNED).

\begin{conjecture}[LNED Conjecture; \cite{zhaogen}]
Let $K$ be a field of characteristic 0, $\mathcal A$ be a $K$-algebra, and $\delta$ be a locally nilpotent $K$-derivation or $K$-$\varepsilon$-derivation of $\mathcal A$.  Then for every ideal $I$ of $\mathcal A$, the image $\delta(I)$ of $I$ under $\delta$ is a Mathieu-Zhao subspace of $\mathcal A$.
\end{conjecture}

We will spend the next section defining and exploring Mathieu-Zhao subspaces for commutative associative algebras, but it should be clear at this point that the Mathieu-Zhao subspace is a potentially very powerful structure.  However, relatively little is still known about this structure even within associative algebras, and this paper marks one of the first formal attempts to study their nature beyond the realm of associative algebras.

\begin{center}
\textbf{Generalizing Ideals}
\end{center}

A Mathieu-Zhao subspace of a ring-algebra is a generalization of an ideal which can be described as a metaphorical "black hole."  Where the elements of an ideal "suck" all other elements straight into the ideal, never to escape again, elements of a Mathieu-Zhao subspace are more like indicators of an "event horizon," meaning that we can only get so close before we are "sucked in" anyway.  We state the formal definition of a Mathieu-Zhao subspace (sometimes simply called a "Mathieu subspace" in the literature) below.  In the tradition of powerful mathematical tools, this definition is deceptively simple, only requiring knowledge of multiplication on a given algebra to begin study.
\theoremstyle{definition}
\begin{definition}[Mathieu-Zhao Subspace; \cite{zhaogen}]
Let $R$ be a commutative unital ring, and $\mathcal{A}$ be a commutative $R$-algebra.  We call an $R$-subspace $M$ a \textit{Mathieu-Zhao subspace} (MZ subspace) if for any $a,b\in\mathcal{A}$ with $a^m\in M$ for any $m\geq 1$, we have $ba^m\in V$ when $m\gg 0$.  That is, there exists $N\geq 1$ (which, in general, depends on $a$ and $b$) such that $ba^m\in M$ for any $m\geq N$.
\end{definition}

We have restricted this definition to \textit{commutative} (associative) algebras for the sake of brevity.  It is easy to generalize this definition to non-commutative cases as have been studied in (\cite{zhaoidempotent}), for example.  We will see in the case of vertex algebras that this question of commutativity requires careful consideration.

We also note that this structure is just a generalization of an ideal.  If we declare $m=1$ in our definition, we have exactly the definition of an ideal.  This is noteworthy because we have countless examples of Mathieu-Zhao subspaces, along with immeasurable theory on those examples, for free.  Let us now observe some examples of non-ideal MZ subspaces.
\\
\begin{example}[\cite{zhaotalk}]
Viewing $\mathbb Q$ as a $\mathbb Z-$algebra, the set
$$M_c:=\Big\{\frac{ac}{b}\Big|a,b\in\mathbb Z; b\neq 0; (b,c)=1\Big\}$$
is an MZ subspace of $\mathbb Q$.
\end{example}

\begin{example}[\cite{fpyz,pakovich}]
The set
$$M:=\Big\{f(x)\in\mathbb C[x]\Big|\int_0^1 f(x)dx=0\Big\}$$
is an MZ subspace of $\mathbb{C}[x]$.
\end{example}

\begin{example}[\cite{zhaogen}]
Let $A=\mathbb C[t,t^{-1}]$ be the ring of Laurent polynomials over $\mathbb C$, and for $\lambda\in\mathbb C$, let $D_{\lambda}=\partial_t+\lambda t^{-1}$.  Then $$M_\lambda:=D_\lambda(A)$$ is an MZ subspace if and only if $\lambda\notin \mathbb Z$ or $\lambda=-1$. 
\end{example}

These examples are all highly nontrivial and demonstrate that, while sometimes hard to pin down, MZ subspaces are everywhere.  In fact, we may informally recognize that there are at least several orders of magnitude more MZ subspaces than there are ideals.
\\
\begin{definition}[\cite{zhaoidempotent}]
Let $\mathcal{A}$ be a (commutative) $R$-algebra and $M$ be a subspace of $\mathcal{A}$.  We define the \textit{radical of $M$} to be the set
$$\mathfrak{r}(M)=\{a\in \mathcal{A}|a^m\in M \text{ for all } m\gg 0\}.$$
In the literature, this set is sometimes denoted by $\sqrt{M}$.

We also define the \textit{strong radical of $M$} to be the set
$$\mathfrak{sr}(M)=\{a\in \mathcal{A}|\text{ for all }b\in \mathcal{A}, ba^m\in M\text{ when } m\gg 0\}.$$
Note that our minimal $m$ may depend on both $a$ and $b$.
\end{definition}

Moving toward how these sets directly connect to MZ subspaces, we have the following immediate consequence:
\\
\begin{lemma}[\cite{intromathieu}]
For any subspace $M$ of a (unital) $R$-algebra $\mathcal{A}$, we have
$$\{0\}\subseteq\mathfrak{sr}(M)\subseteq \mathfrak{r}(M).$$

\begin{proof} We have $0\in\mathfrak{sr}(M)$ from the definition of strong radical, and that $M$ is a subspace (thus $0\in M$).

Now, let $a\in\mathfrak{sr}(M)$ be given.  Then for any $b\in \mathcal A$, there exists $N\in\mathbb{N}$ such that for all $m\geq N$, we have
$$ba^m\in M.$$
Since $\mathcal A$ is unital, take $b=1_\mathcal{A}$.  Then we have
$$1_\mathcal{A}\cdot a^m=a^m\in M$$
for all sufficiently large $m$. \end{proof}
\end{lemma}

The notion of the radical and the strong radical will become especially important in the context of vertex algebras.  One simple reason for this importance is that we intend to extend the following theorem:
\\
\begin{theorem}[\cite{intromathieu}]
Let $\mathcal A$ be an $R$-algebra and $M$ be a subspace of $\mathcal A$.  Then $M$ is a Mathieu-Zhao subspace of $\mathcal A$ if and only if $$\mathfrak{r}(M)=\mathfrak{sr}(M).$$

\begin{proof} $(\Rightarrow)$ Let $M$ be an MZ subspace of $\mathcal A$.  By Lemma 10, we need only to show that $$\mathfrak{r}(M)\subseteq \mathfrak{sr}(M).$$
Let $a\in\mathfrak{r}(M)$ be given; i.e., there exists $N_1\in\mathbb N$ such that $a^m\in M$ for all $m\geq N_1$.  This implies that $(a^{N_1})^t\in M$ for all $t\geq 1$.  Since $M$ is MZ, we know that for all $b\in\mathcal A$, there exists $N\in\mathbb N$ such that $b(a^{N_1})^t\in M$ for all $t\geq N$.  We can easily see that we may repeat this argument for $N_i:=N_1+i$ for all $0\leq i<N_1$ since $a\in\mathfrak r (M)$ to cover all sufficiently large integers $m$.  In more concrete terms, this implies that for all $b\in\mathcal A$, we have $ba^m\in V$ when $m\gg 0$ which implies that $a\in\mathfrak{sr}(M)$.

$(\Leftarrow)$ Now suppose that $\mathfrak{r}(M)= \mathfrak{sr}(M)$ (which, again, is equivalent to supposing that $\mathfrak{r}(M)\subseteq \mathfrak{sr}(M)$).  Let $a\in\mathcal A$ be such that $a^m\in M$ for all $m\geq 1$.  Clearly, we must have $a\in\mathfrak{r}(M)=\mathfrak{sr}(M)$, which implies that for any $b\in \mathcal A$, there exists $N\in\mathbb{N}$ such that $ba^m\in V$ for all $m\geq N$.  Then, by definition, $M$ is an MZ subspace of $\mathcal A$. \end{proof}
\end{theorem}

This is a powerful theorem which leads, through the following corollaries, to a nearly trivial classification for a large set of MZ subspaces which may not be ideals (\cite{intromathieu}).  For the following, let $\mathcal A$ be an $R$-algebra and $M$ be a subspace of $\mathcal A$.
\\
\begin{corollary}
If $\mathfrak{r}(M)=\{0\}$, then $M$ is an MZ subspace of $\mathcal A$.

\begin{proof} We have $$\{0\}\subseteq \mathfrak{sr}(M)\subseteq\mathfrak{r}(M)=\{0\}.$$ \end{proof}
\end{corollary}

\begin{corollary}
If $\mathfrak{sr}(M)=\mathcal A$, then $M$ is an MZ subspace of $\mathcal A$.

\begin{proof} We have the similar inclusion $$\mathcal A=\mathfrak{sr}(M)\subseteq \mathfrak{r}(M)\subseteq \mathcal A.$$\end{proof}
\end{corollary}

\begin{corollary}
If $1\in\mathfrak{sr}(M)$, then $M=\mathcal A$.

\begin{proof} We know that for any $m\in\mathbb N$, any $b\in \mathcal A$ yields $$b\cdot 1^m=b\in V.$$\end{proof}
\end{corollary}

Before closing this section, we provide one theorem, which will become vital in our final chapter.  This theorem concerns the specific case of MZ subspaces of the polynomial algebra $K[x]$ over a field $K$ of characteristic zero and can be generalized to the algebras $K[x_i|0\leq i\leq n]$ and $K[x_i|i\in\mathbb N]$.
\\
\begin{theorem}[\cite{zhaoimages}]
Let $K$ be a field of characteristic zero, $\{n_i|i\geq 1\}$ be a strictly increasing sequence such that $0\neq n_i\in\mathbb N$ for all $i\geq 1$ and $n_{i+1}-n_i\neq 1$ for infinitely many $i\geq 1$, and $M$ be the $K$-subspace of $K[x]$ spanned by $\{x^{n_i}|i\geq 1\}$ over $K$.  Then the following are equivalent:
\begin{itemize}
\item[(1)] $\mathfrak{r}(M)=\{0\}$;
\item[(2)] $M$ is a Mathieu-Zhao subspace of $K[x]$;
\item[(3)] there exists no $d\in\mathbb{N}$ such that $md\in\{n_i|i\geq 1\}$ for all $m\geq 1$.
\end{itemize}

\begin{proof} $(1)\Rightarrow (2)$: See Corollary 12.
\\
$(2)\Rightarrow (3)$: Assume by contradiction that there exists $d\geq 1$ such that $md\in\{n_i|i\geq 1\}$ for all $m\geq 1$.  This implies that $x^{md}\in M$ for all $m\geq 1$.  We know that $d>1$; otherwise, we would have $\{n_i|i\geq 1\}=\mathbb N$, a contradiction.  Furthermore, since $(x^d)^m=x^{dm}\in M$ for all $m\geq 1$ and since $M$ is MZ in $K[x]$, we know that for each $0\leq r\leq d-1$, there exists $N_r\geq 1$ such that for all $m\geq N_r$, we have $x^{md+r}=(x^d)^mx^r\in M$ by the division algorithm.  Allowing $N=\text{max}\{N_r|0\leq r\leq d-1\}$, we see that for all $k\geq Nd$, we have $x^k\in M$ which implies that $k\in\{n_i|i\geq 1\}$, contradicting our assumption on $\{n_i|i\geq 1\}$ (namely, that $n_{i+1}-n_i\neq 1$ for infinitely many $i$).  We conclude that no such $d$ exists, thus $(2)\Rightarrow (3)$.
\\
$(3)\Rightarrow (1)$: Assume by contradiction that there exists $0\neq f(x)\in\mathfrak{r}(M)$.  That is, $f^m(x)\in M$ when $m\gg 0$.  Since $\{n_i|i\geq 1\}$ consists of strictly positive integers, we know that $1\notin M$.  Thus $c:=\text{deg}(f(x))\geq 1$, and furthermore $x^{cm}\in M$ when $m\gg 0$ by definition of $M$ (as a space spanned by monomials).  But this means that there exists $N\in\mathbb N$ such that $mc\in\{n_i|i\geq 1\}$ for all $m\geq N$.  Then considering $d:=N\cdot c$, we have $md\in\{n_i|i\geq 1\}$ for all $m\geq 1$, a contradiction.  We conclude that no such $f$ exists; thus, $(3)\Rightarrow (1)$.
\end{proof}
\end{theorem}

\begin{example}
Consider the polynomial ring $\mathbb C [x]$ as a $\mathbb C$-algebra.  We fix $k>1$ and introduce a map
$$\phi:\mathbb{C}[x]\to\mathbb{C}[x],$$
$$x^n\mapsto e^{\frac{2\pi in}{k}}x^n.$$
We may alternatively view this map as the group action on $\mathbb C[x]$ by $\mathbb Z/k\mathbb Z$, the cyclic group of $k$ elements, and thus partition $\mathbb{C}[x]$ into eigenspaces
$$\mathbb{C}[x]=\bigoplus_{l=0}^{k-1}\mathbb{C}[x]^l,$$
where
\[
\begin{array}{l c l}
	\mathbb{C}[x]^0 & = & \{f(x)\in\mathbb{C}[x]|f(x)=\sum_{m}c_mx^{km}\text{ where }c_m\in\mathbb C\}, \\
	\mathbb{C}[x]^1 & = & \{f(x)\in\mathbb{C}[x]|f(x)=\sum_{m}c_mx^{km+1}\text{ where }c_m\in\mathbb{C}\}, \\
	& \vdots & \\
	\mathbb{C}[x]^{k-1} & = & \{f(x)\in\mathbb{C}[x]|f(x)=\sum_{m}c_mx^{km+(k-1)}\text{ where }c_m\in\mathbb{C}\}.
\end{array}
\]
In a geometric sense, $\mathbb{C}[x]^0$ is the set of polynomials fixed by our map $\phi$ (i.e. are rotated by 0), $\mathbb{C}[x]^1$ is the set of polynomials rotated by $\frac{2\pi i}{k}$, $\mathbb{C}[x]^2$ is the set of polynomials rotated by $\frac{4\pi i}{k}$, and so on.
\end{example}

\begin{proposition}
Fix $k>1$, and consider the decomposition
$$\mathbb C[x]=\bigoplus_{l=0}^{k-1}\mathbb C[x]^l.$$ 
Then $\mathbb{C}[x]^0$ is \textit{not} a Mathieu-Zhao subspace of $\mathbb{C}[x]$.

\begin{proof} Assume $\mathbb{C}[x]^{(0)}$ is a Mathieu-Zhao subspace of $\mathbb{C}[x]$.  Note that 1 is a fixed point of our map $\phi$ (regardless of our choice for $k$).  This implies that $1\in\mathbb C[x]^{(0)}$.  It naturally follows that $1\in\mathfrak{r}(\mathbb C[x]^0)$ since $1^n=1$ for all $n\in\mathbb Z$.  As $\mathbb{C}[x]^0$ is MZ in $\mathbb{C}[x]$, we know that $1\in\mathfrak{sr}(\mathbb{C}[x]^0)$, but by Corollary 14, this implies that $\mathbb{C}[x]^0=\mathbb{C}[x]$, a contradiction.  Hence, $\mathbb{C}[x]^0$ must not be MZ in $\mathbb{C}[x]$.\end{proof}
\end{proposition}

\begin{proposition}
Again, fix $k>1$ in the example above.  Then the subspace $$M:=\bigoplus_{l=1}^{k-1}\mathbb{C}[x]^{l_r}$$ is a Mathieu-Zhao subspace of $\mathbb{C}[x]$.

\begin{proof}As in the proof of Proposition 17, we first note that $M$ is spanned by the set of monomials $\{x^{l+km}|m\geq 1\}$.  Again, we naturally associate a sequence $\{l+km|m\geq 1\}\{1,2,\dots,k-1,k+1,\dots\}$ which satisfies our condition in Theorem 15.  Note that this sequence lacks all multiples of $k$ (thus the difference of consecutive terms $km+1-(km-1)=2\neq 1$ for all $m\in\mathbb N$).

Having satisfied condition (3) of Theorem 15, we conclude $M$ is MZ in $\mathbb C[x]$.\end{proof}
\end{proposition}

In moving toward vertex algebras, this theorem and the subsequent example should give us an intuition about what we will see.  Our aim is to exploit natural connections between associative algebras and vertex algebras, which, as we will see, are just generalizations of those same associative algebras, albeit with more technicalities to consider.

\newpage

\begin{center}
CHAPTER II: VERTEX ALGEBRAS AND THEIR MATHIEU-ZHAO SUBSPACES
\end{center}

\begin{center}
\textbf{Vertex Algebras}
\end{center}

In the following, we formally define a vertex algebra (VA).  Informally, vertex algebras generalize associative algebras in ways that allow for the novel types of dimensional analysis employed in (\cite{moonshine}).

\begin{definition}[Vertex Algebra; \cite{lepowskyli, borcherds1, flm}]
Let $z,z_0,z_1,z_2$ be indeterminates, and $V=\coprod_{n\in\mathbb{Z}}V_n$ a $\mathbb{Z}$-graded $\mathbb C$-vector space, wherein for $v\in V_n$, we have $n=\text{wt }v$.  Equip $V$ with an injective map
\[
\begin{array}{l l l l}
Y:&V &\to& (\text{End }V)[[z,z^{-1}]]\\
&v&\mapsto &Y(v,z)=\sum_{n\in\mathbb{Z}}v(n)z^{-n-1}
\end{array}
\]
and with a vacuum element $\textbf{1}\in V_0$ satisfying the following:
\[
\begin{array}{l}
u(n)v=0,\text{ for sufficiently large }n;\\
Y(\textbf{1},z)=1;\\
Y(v,z)\textbf{1}\in V[[z]],\text{ and } \lim_{z\to 0}Y(v,z)\textbf{1}=v;\text{ In particular, }v(-1)\textbf{1}=v;\\
z_0^{-1}\delta\Big(\frac{z_1-z_2}{z_0}\Big)Y(u,z_1)Y(v,z_2)-z_0^{-1}\delta\Big(\frac{z_2-z_1}{z_0}\Big)Y(v,z_2)Y(u,z_1)\\
\hspace{1.5cm}=z_2^{-1}\delta\Big(\frac{z_1-z_0}{z_2}\Big)Y(Y(u,z_0)v,z_2),
\end{array}
\]
where
$$\delta(z)=\sum_{n\in\mathbb Z}z^n.$$
This last axiom is called the Jacobi identity for vertex algebras.  Then $V$ is a \textit{vertex algebra}, denoted by $(V,Y,\textbf{1})$, although we may simply denote it by $V$ when such a labeling is unambiguous.
\end{definition}

For the primary aims of this paper, we only require the above axioms.  However, it may be easier, in certain situations, to distinguish a vector $\omega\in V$, which classifies a specific type of vertex algebra called a \textit{vertex operator algebra} (VOA), defined as follows:

\begin{definition}[Vertex Operator Algebra; \cite{lepowskyli, borcherds1, flm}]
Let $V=\coprod_{n\in\mathbb{Z}}V_n$ be a $\mathbb{Z}$-graded $\mathbb C$-vector space, wherein for $v\in V_n$, we have $n=\text{wt }v$ such that $\text{dim }V_n<\infty$.  Suppose there exists a vector $\omega\in V$ satisfying the following:
\[
[L(m),L(n)]= (m-n)L(m+n)+\frac{1}{12}\delta_{m+n,0}(m^3-m)\delta_{m+n,0}(\text{rank }V)
\]
for $m,n\in\mathbb Z$, where
$$L(n)=\omega_{n+1}, \text{ for }n\in\mathbb Z,\text{ i.e., }Y(\omega,z)=\sum_{n\in\mathbb Z}L(n)z^{-n-2}$$
and
\[ \begin{array}{l}
\text{rank }V\in\mathbb Q;\\
L(0)v=nv=(\text{wt }v)v\text{ for }v\in V_n;\\
\frac{d}{dz}Y(v,z)=Y(L(-1)v,z).
\end{array}
\]
Then $V$ is a \textit{vertex operator algebra}, denoted by $(V,Y,\textbf{1},\omega)$, although, as before, we may simply denote it by $V$ when such a labeling is unambiguous.
\end{definition}

\begin{definition}
Let $V$ be a vertex algebra.  An \textit{ideal} $I$ of $V$ is a subspace of $V$ such that $v(n)u\in I$ and $u(n)v\in I$ for any $u\in V,v\in I,n\in\mathbb Z$.
\end{definition}

\begin{definition}
Let $V$ be a vertex algebra.  $V$ is said to be a \textit{simple vertex algebra} if $V$ contains no nontrivial strict ideals.
\end{definition}

\begin{example}
Consider the polynomial ring $\mathbb C[t,t^{-1}]$ equipped with a differential $\partial$.  Define for $f,g\in\mathbb C[t,t^{-1}]$,
$$Y(f,z)g=(e^{z\partial}f)g,$$
i.e.,
$$f(-n-1)g=\frac{\partial^n(f)}{n!} g.$$
Then $(\mathbb C[t,t^{-1}],Y,1)$ is a unital commutative vertex algebra.\qed
\end{example}

In fact, any associative algebra with a well-defined differential is a vertex algebra.

\begin{definition}
We define a \textit{Lie algebra} to be a vector space, \textbf{h}, over $\mathbb C$, together with a bilinear operation
\[\begin{array}{l l l}
\textbf{h}\times\textbf{h}&\to& \textbf{h}\\
(x,y)&\mapsto &[x,y],
\end{array}
\]
such that for all $a,b,c\in\textbf h$
$$[a,[b,c]]+[c,[a,b]]+[b,[a,c]]=0$$
and
$$[a,b]=-[b,a].$$

Furthermore, let $\textbf{h}$ be a finite-dimensional Lie algebra over $\mathbb C$, with basis $\{x_1,\dots,x_m\}$, and bilinear operation
$$[x_i,x_j]=\sum_{k=1}^mc_{i,j,k}x_k,$$
where all $c_{i,j,k}$ are elements of $K$.  Then the \textit{universal enveloping algebra} of \textbf{h}, denoted $U(\textbf{h})$, is the associative algebra generated by $\{x_1,\dots,x_m\}$ subject exclusively to the relations
$$[x_i,x_j]=x_ix_j-x_jx_i=\sum_{k=1}^lc_{i,j,k}x_k.$$
\end{definition}

\begin{definition}[\cite{dongintro, flm}]
Let $\textbf{h}$ be a vector space equipped with a symmetric nondegenerate bilinear form $\langle\cdot,\cdot\rangle$.  Viewing \textbf{h} as an abelian Lie algebra, consider the affine Lie algebra
$$\hat{\textbf{h}}:=\textbf{h}\otimes\mathbb C[t,t^{-1}]\oplus\mathbb Cc$$
(here, we call $c$ the \textit{central element}) with
$$[x\otimes t^m, y\otimes t^n]=\langle x,y\rangle m\delta_{m+n,0}c$$
for $x,y\in\textbf{h}, m,n\in\mathbb Z$, and
$$[c,\hat{\textbf{h}}]=0.$$
Here, $\delta_{x,y}$ is the Kronecker $\delta$ function, i.e.,
\[\delta_{a,b}=\begin{cases}
	1&a=b\\
    0&a\neq b.
\end{cases}
\]
Then set
$$\hat{\textbf{h}}^+=\textbf{h}\otimes t\mathbb C[t], \hat{\textbf{h}}^-=\textbf{h}\otimes t^{-1}\mathbb C[t^{-1}].$$
Then we call the subalgebra
$$\hat{\textbf{h}}_{\mathbb Z}:=\hat{\textbf{h}}^+\oplus\hat{\textbf{h}}^-\oplus\mathbb Cc$$
a \textit{Heisenberg algebra}.  Consider the induced $\hat{\textbf{h}}$-module
$$M(1):=U(\hat{\textbf{h}})\otimes_{U(\textbf{h}\otimes\mathbb C[t]\oplus\mathbb Cc)}\mathbb C_0$$
with $\textbf{h}\otimes t\mathbb c[t]$ acting trivially on $\mathbb C$, $\textbf{h}$ acting as $\langle \textbf{h},0\rangle$, and $c$ acting as 1.  Now, for $\alpha\in\textbf{h}$ and $n\in\mathbb Z$, write $\alpha(n)=\alpha\otimes t^n$, set
$$\alpha(z)=\sum_{n\in\mathbb Z}\alpha(n)z^{-n-1},$$
and define
$$Y(v,z)=\leftidx{_\circ^\circ}{\Big(}\frac{1}{(n_1-1)!}\Big(\frac{d}{dz}\Big)^{n_1-1}\alpha_{i_1}(z)\Big)\cdots\Big(\frac{1}{(n_k-1)!}\Big(\frac{d}{dz}\Big)^{n_k-1}\alpha_{i_k}(z)\Big)_\circ^\circ,$$
where $$v=\alpha_{1}(-n_1)\cdots\alpha_{k}(-n_k)\textbf{1}$$
for $\alpha_1,\dots,\alpha_k\in\textbf h, n_1,\dots,n_k\in\mathbb N$, and open colons indicate a \textit{normal ordering}, i.e., we may reorder our operators such that all of $\alpha(n) (\alpha\in\textbf{h},n<0)$ are to the left of all $\alpha(n) (\alpha\in\textbf{h}, n\geq 0)$ before evaluation.  We extend this definition to all $v\in M(1)$ linearly.  Set $\textbf{1}=1,\omega=\frac{1}{2}\sum_{i=1}^d\alpha_i(-1)^2\in M(1)$ where $\{\alpha_i|1\leq i\leq d\}$ is an orthonormal basis for \textbf{h}.  Then we call $(M(1),Y,\textbf{1},\omega)$ the \textit{Heisenberg vertex operator algebra}.  In particular, if $|\{\alpha_i|1\leq i\leq d\}|=d$, we denote the corresponding Heisenberg VOA by $M_d(1)$.
\end{definition}

For the remainder of this paper, we will focus on the one-dimensional Heisenberg VOA, $M_1(1)$, and will always denote the basis of $\textbf{h}$ by $\{\alpha\}$.  In a sense, our effort here is to begin to transfer what we know of MZ subspaces of the polynomial algebra $\mathbb C[x]$ to what can loosely be described as its partner VOA, $M_1(1)$.  When we formally define MZ subspaces for vertex algebras and begin to develop their theory, we will consequently solidify this "partnership" between $\mathbb C[x]$ and $M_1(1)$.

\begin{theorem}[\cite{flm}]
$M(1)$ is a simple vertex algebra.
\end{theorem}

The following lemma, stated as a simple matter of fact in (\cite{moonshine}) and proven in detail in (\cite{lepowskyli}), gives a property we will call \textit{pseudo-associativity}.  The equation therein is derived directly from the Jacobi Identity.

\begin{lemma}[Pseudo-associativity; \cite{borcherds1,lepowskyli}] Consider $M_1(1)$.  For any vectors $u,v,w\in M_1(1)$ and any integers $m,n\in\mathbb{Z}$, we have
$$\big(u(m)v\big)(n)w=\sum_{i\geq 0}(-1)^i\binom{m}{i}\big(u(m-i)v(n+i)w-(-1)^mv(m+n-i)u(i)w\big).$$
\end{lemma}

\begin{definition} Consider an element $A\in M_1(1)$ such that
$$A=\alpha(-n_1)\cdots\alpha(-n_d)\textbf{1}$$
for some $n_1,\dots,n_d\in\mathbb N$.  Then we say that $d$ is the \textit{$\alpha$-length of $A$}.  Moreover, we denote
$$M_1^d(1):=\text{span}\{\alpha(-n_1)\cdots\alpha(-n_d)\textbf{1}\text{ }|\text{ }n_1,\dots,n_d\in\mathbb N\}.$$
\end{definition}

\begin{lemma}[Pseudo-commutativity]
Consider $M_1(1)$, and let $w\in M_1(1),n,m\in\mathbb Z$.  Then
$$\alpha(n)\alpha(m)w=\alpha(m)\alpha(n)w+n\langle \alpha,\alpha\rangle \delta_{m+n,0}w,$$
where $\delta_{x,y}$ again denotes the Kronecker $\delta$ function.
\end{lemma}

\begin{lemma} Let $$A:=\alpha(-1)\cdots\alpha(-1)\textbf{1}\in M_1^d(1)$$ be given.  Then we have $$A(q_1)w\in \bigoplus_{i=0}^{k+d}M_1^i(1)$$ for all $q_1\in\mathbb Z$ and $w\in M_1^k(1)$.  In particular, we have
$$A(q_1)w\in\bigoplus_{j=0}^dM^{k+d-2j}.$$

\begin{proof} Note first by Lemma 29 that for $n\geq 1$, we have
$$\alpha(n)\underbrace{\alpha(-n)\cdots\alpha(-n)}_{p\text{ times}}\textbf{1}=pn\underbrace{\alpha(-n)\cdots\alpha(-n)}_{p-1\text{ times}}\textbf{1}.$$
We proceed by induction on $d$.  First, let $d=1$.  Then we have
$$A(q_1)w=\big(\alpha(-1)\textbf{1})(q_1)w=\alpha(q_1)w\in M_1^{k-1}(1)\oplus M_1^{k+1}(1)\subseteq \bigoplus_{i=0}^{k+1}M_1^i(1).$$
Now, let $d=2$.  We have
$$A(q_1)w=\big(\alpha(-1)\alpha(-1)\textbf{1})(q_1)w.$$
To apply pseudo-associativity, set $u=\alpha(-1)\textbf{1}, v=\alpha(-1)\textbf{1}, w=w, m=-1,n=q_1$.  Then we have
\[
\begin{array}{l l l}
A(q_1)\textbf{1}&=&\sum_{i\geq 0}(-1)^i\binom{-1}{i}\Big(\alpha(-1-i)\big(\alpha(-1)\textbf{1}\big)(q_1+i)w\\
&&-(-1)^{-1}\big(\alpha(-1)\textbf{1}\big)(-1+q_1-i)\alpha(i)w\Big).
\end{array}
\]
Then, we apply the same logic as we did in the case of $d=1$ and note that our component parts have $\alpha$-length either $k+2$, $k$, or $k-2$.  In other words, we have
$$A(q_1)w\in M_1^{k-2}(1)\oplus M_1^k(1)\oplus M_1^{k+2}(1)\subseteq \bigoplus_{i=0}^{k+2}M_1^i(1).$$

Now suppose that our assumption holds for $d=s$, and consider $d=s+1$.  Set
$$B:=\underbrace{\alpha(-1)\cdots\alpha(-1)}_{s\text{ times}}\textbf{1}.$$
By our induction hypothesis, we have
\[\begin{array}{l l l}
\big(\underbrace{\alpha(-1)\cdots\alpha(-1)}_{s+1\text{ times}}\textbf{1}\big)(q_1)w&=&\sum_{i\geq 0}\big(\alpha(-1)\textbf{1}\big)(-1-i)B(q_1+i)w\\
&&+B(-1+q_1+i)\big(\alpha(-1)\textbf{1}\big)(i)w\\
&\in & \bigoplus_{i=0}^{k+s+1}M_1^i(1).
\end{array}\]
\end{proof}
\end{lemma}

\begin{lemma} Let $p,d\in\mathbb{N}$ be such that $pd\geq 2$, and let $A\in M_1(1)$ be of the form
$$A=\underbrace{\alpha(-1)\cdots\alpha(-1)}_{d\text{ times}}\textbf{1}.$$
Then for $w\in M_1^k(1)$, we have
$$A(q_1)\cdots A(q_p)w\in\bigoplus_{i=0}^{pd}M_1^{pd+k-2i}(1),$$
for all $q_1,\dots,q_p\in\mathbb Z$.

\begin{proof} We proceed by induction on $p$.  The case of $p=1$ is exactly Lemma 30.

Now, assume
$$B:=A(q_1)\cdots A(q_p)w\in\bigoplus_{i=0}^{pd}M_1^{pd+k-2i}(1).$$
Then we consider the element
\[\begin{array}{l l l}
A(q_0)\cdots A(q_p)w&=&\big(\underbrace{\alpha(-1)\cdots\alpha(-1)}_{d\text{ times}}\textbf{1}\big)(q_0)B\\
&=& \sum_{j\geq 0}\alpha(-1-i)\big(\underbrace{\alpha(-1)\cdots\alpha(-1)}_{d-1\text{ times}}\textbf{1}\big)(q_0+i)B\\
&& +\big(\underbrace{\alpha(-1)\cdots\alpha(-1)}_{d-1\text{ times}}\textbf{1}\big)(q_0-1-i)\alpha(i)B.
\end{array}
\]
We are now able to apply Lemma 30 component-wise and see that, by our induction hypothesis, each of the terms in this equation constitutes an element of $\bigoplus_{i=0}^{pd}M_1^{pd+k-2i}(1)$.
\end{proof}
\end{lemma}

\begin{center}
\textbf{Main Results}
\end{center}

\addcontentsline{toc}{subsection}{Main Results}
Recall Definition 9.  Therein, we defined the concepts of radical and strong radical for an associative algebra.  We will first define these concepts for vertex algebras and then use these concepts to anchor our understanding of MZ subspaces within VA's.  Unlike in the case of associative algebras, Mathieu-Zhao subspaces of vertex algebras do not lend themselves to a concise and easily understood definition without such sets.

\begin{definition}
Let $(V,Y,\textbf{1})$ be a vertex algebra and $M$ be a subspace of $V$.  Then define the \textit{radical of $M$} to be
\[
\begin{array}{l l l}
	\mathfrak{r}(M)&=&\{v\in V|\text{ there exists }m\geq 1,\text{ such that }v(n_1)\cdots v(n_t)\textbf{1}\in M, \\
	& & \hphantom{space}\text{ for all }t\geq m,\text{ with }n_1,\dots,n_t\in\mathbb Z\}.
\end{array}
\]
Similarly, we define the \textit{left strong radical of $M$} to be
\[
\begin{array}{l l l}
	l\mathfrak{sr}(M)&=&\{v\in V|\text{ there exists }m\geq 1,\text{ such that } b(s) v(n_1)\cdots v(n_t)\textbf{1}\in M \\
	&&\hphantom{space}\text{ for all }t\geq m,b\in V,s,n_1,\dots,n_t\in\mathbb Z\},
\end{array}
\]
and define the \textit{right strong radical of $M$} to be
\[
\begin{array}{l l l}
	r\mathfrak{sr}(M)&=& \{v\in V|\text{ there exists } m\geq 1,\text{ such that }\big(v(n_1)\cdots v(n_t)\textbf{1}\big)(n)w\in M,\\
	&&\hphantom{space}\text{ for all } t\geq m,w\in V,n,n_1,\dots,n_t\in \mathbb Z\}.
\end{array}
\]
Then define the \textit{strong radical of $M$} to be
$$\mathfrak{sr}(M)=l\mathfrak{sr}(M)\cap r\mathfrak{sr}(M).$$
We say that $M$ is a \textit{Mathieu-Zhao subspace of $V$} if
$$\mathfrak{r}(M)=\mathfrak{sr}(M).$$
\end{definition}

In the interest of simplifying notation, the following tool will allow us to forego the qualifiers \textit{left} and \textit{right} strong radicals in certain circumstances (see Proposition 34).

\begin{lemma}[Skew Symmetry]
Let $D$ denote the derivation on $M_1(1)$ defined by
$$D(v)=v(-2)\textbf{1}.$$
By the definition of a vertex algebra, we have
$$Y(B,z)A=e^{zD}Y(A,-z)B.$$
Then we observe that
\[
\begin{array}{r c l}
	e^{zD}Y(A,-z)B & = & \sum_{m\in\mathbb Z}\sum_{i=0}^\infty\frac{z^iD^i}{i!}A(m)(-z)^{-m-1}B \\
	&& \\
	& = & \sum_{n\in\mathbb Z}\sum_{i=0}^\infty(-1)^{-n-i-1}\frac{D^i}{i!}A(n+i)Bz^{-n-1} \\
	&&\\
	\Rightarrow B(n)A & = & \sum_{i=0}^\infty (-1)^{-n-i-1}\frac{D^i}{i!}A(n+i)B.
\end{array}
\]
\end{lemma}

\begin{proposition}
Let $M$ be a subspace of $V$ such that $D(M)\subseteq M$.  Then
$$\mathfrak{sr}(M)=l\mathfrak{sr}(M)=r\mathfrak{sr}(M).$$

\begin{proof}  First, let $v\in l\mathfrak{sr}(M)$, and set $B=v(n_1)\cdots v(n_t)\textbf{1},A=w$.  Then we have
$$A(m)B=w(m)v(n_1)\cdots v(n_t)\textbf{1}\in M.$$
By Lemma 29 and the fact that $D(M)\subseteq M$, this implies that for all $n\in\mathbb Z$, we have
$$B(n)A=\big(v(n_1)\cdots v(n_t)\textbf{1}\big)(n)w\in M\Rightarrow v\in r\mathfrak{sr}(M).$$

Conversely, let $v\in r\mathfrak{sr}(M)$ and set $A=v(n_1)\cdots v(n_t)\textbf{1}, B=b$.  Then we have
$$A(m)B=\big(v(n_1)\cdots v(n_t)\textbf{1}\big)(m)w\in M,$$
which, again, implies that for all $n\in\mathbb Z$, we have the containment
$$w(n)v(n_1)\cdots v(n_t)\textbf{1}\in M,$$
thus $v\in l\mathfrak{sr}(M)$.\end{proof}
\end{proposition}

We can now state the generalizations of the basic properties of MZ subspaces in our new context.

\begin{lemma}
For any subspace $M$ of a vertex algebra $V$, we have
$$\{0\}\subseteq \mathfrak{sr}(M)\subseteq \mathfrak{r}(M).$$

\begin{proof} Similarly to Lemma 10, since $M$ is a subspace, $0\in\mathfrak{sr}(M)$ is obvious.

Now consider $v\in\mathfrak{sr}(M)$.  There exists $m\geq 1$ such that $$b(s)v(n_1)\cdots v(n_t)\textbf{1}\in M$$ for all $t\geq m,b,w\in V, s,n,n_1,\dots,n_t\in\mathbb Z$.  In particular, consider $b=\textbf{1},s=-1$.  Then observe $$\textbf{1}(-1)v(n_1)\cdots v(n_t)=v(n_1)\cdots v(n_t)\textbf{1}\in M,$$ which implies $v\in\mathfrak{r}(M)$.\end{proof}
\end{lemma}

Just as in the associative case, we are able to immediately observe the following corollaries, whose proofs follow the exact same lines of reasoning as Corollaries 12, 13, and 14:

\begin{corollary}
If $\mathfrak{r}(M)=\{0\}$, then $M$ is an MZ subspace of $V$.\qed
\end{corollary}

\begin{corollary}
If $\mathfrak{sr}(M)=V$, then $M$ is an MZ subspace of $V$.\qed
\end{corollary}

\begin{corollary}
If $\textbf{1}\in\mathfrak{sr}(M)$, then $M=V$.\qed
\end{corollary}

\begin{lemma}
Let $I$ be an ideal in $V$.  Then $I$ is a Mathieu-Zhao subspace of $V$.

\begin{proof}Let $u\in\mathfrak r(I)$.  Then there exists $m\geq 1$ such that $u(n_1)\cdots u(n_t)\textbf{1}\in I$ for all $t\geq m, n_1,\dots,n_t\in\mathbb Z$.  Since $I$ is an ideal, this means that for all $b\in V,s\in\mathbb Z$, we have $b(s)u(n_1)\cdots u(n_t)\textbf{1}\in I$ and $u\in\mathfrak{sr}(M)$.\end{proof}
\end{lemma}

So far, we have seen that at least a significant number of elementary properties of MZ subspaces of associative algebras naturally extend to those of VA's.  However, if we recall Example 23, we are able to observe a more applied extension of theory:

\begin{proposition}
Recall that for $\mathbb C[t,t^{-1}]$ together with a differential $\partial$ (a space which we denote by $(\mathbb C[t,t^{-1}],\partial)$), we define for $f,g\in\mathbb{C}[t,t^{-1}]$,
$$Y(f,z)g=(e^{z\partial}f)g.$$
If $(M,Y,1)$ is a Mathieu-Zhao subspace of $(\mathbb{C}[t,t^{-1}],Y,1)$, then $(M,\partial)$ is a Mathieu-Zhao subspace of $(\mathbb{C}[t,t^{-1}],\partial)$.

\begin{proof}Since $(M,Y,1)$ is an MZ subspace of $(\mathbb C[t,t^{-1}],Y,1)$, whenever we have an element $f$ such that $f(-n_1-1)\cdots f(-n_p-1)1\in (M,Y,1)$ for all $p\gg 0$, we are also guaranteed that any choice of $g\in\mathbb C[t,t^{-1}]$ and $s\in\mathbb Z$ yields
$$g(s)f(-m_1-1)\cdots f(-m_q-1)1\in (M,Y,1)$$
for all $q\gg 0$ (depending on $g$).  In particular, if we choose $s=-1$, $m_1=\cdots=m_q=0$, and $g$ to be any Laurent polynomial, then we have
$$gf^q\in (M,\partial)$$
for all $q\gg 0$ (depending on $g$).  This is exactly the definition of an MZ subspace for an associative algebra.\end{proof}
\end{proposition}

\begin{definition}
Let $V$ be a vertex algebra over $\mathbb C$ and $\mathcal T\subseteq M$ be any subset of a $V$-submodule $M$ (i.e., $M$ as a subspace of $V$ satisfies the vertex algebra axioms under the same $Y$ and $\textbf{1}$ as $V$).  Then the \textit{annihilating space of $\mathcal T$} is
$$\text{Ann}_V(\mathcal T)=\{v\in V|\text{ for all }u\in\mathcal T, Y(u,z)v=0\}.$$
\end{definition}

For a complete proof of the following lemma, see (\cite{lepowskyli}).  The proof therein requires a property called the weak associativity of VA's, a concept we will otherwise not need in this paper.

\begin{lemma}
Let $\mathcal T\subseteq M$ be a subset of a $V$-submodule $M$.  Then $\text{Ann}_V(\mathcal T)$ is an ideal of $V$.  Moreover, if we denote the submodule generated by the elements of $\mathcal T$ by $\langle\mathcal T\rangle$, then
$$\text{Ann}_V(\mathcal T)=\text{Ann}_V(\langle\mathcal T\rangle).$$\qed
\end{lemma}

\begin{corollary}
If $V$ is a simple vertex algebra, then
$$\text{Ann}_V(\mathcal T)=\{0\}.$$\qed
\end{corollary}

The following example will be similar to Example 16.  All of our main results will be motivated by this example, but we will see that more specific contexts (i.e., restricting our choice of vertex algebra, $V$) allow us to employ different proof strategies.  In the development of a thorough theoretical survey of MZ subspaces of VA's, such an excess of proof methods is certainly desired.

\begin{example}
Consider $M_1(1)$.  We fix $k>1$, and partition our space $M_1(1)$ as follows: Define
\[\begin{array}{l c l c}
g:&M_1(1)&\to&M_1(1)\\
& \underbrace{\alpha(-1)\cdots\alpha(-1)}_{l\text{ times}}\textbf{1}&\mapsto&\Big(\text{exp}\big(\frac{2\pi i}{k}\big)\Big)^l\underbrace{\alpha(-1)\cdots\alpha(-1)}_{l\text{ times}}\textbf{1}.
\end{array}
\]
Then the direct sum
$$M_1(1)=\bigoplus_{l=0}^{k-1}M_1^{(l)}(1),$$
where
\[M_1^{(l)}(1):=\bigoplus_{m\in\mathbb Z}M_1^{l+km}(1),\]
gives an eigenspace decomposition of $M_1(1)$.
\end{example}

\begin{theorem}[Main Result 1]
Let $V$ be a simple vertex algebra, $k>1$ be fixed and $g$ be an automorphism of $V$ of order $k$.  Then we may write
$$V=\bigoplus_{l=0}^{k-1}V^l$$
where
$$V=\{v\in V\text{ }|\text{ }g(v)=e^{\frac{2\pi il}{k}}v\}.$$
For $0\leq i\leq k-1$, the subspace $V^i$ is a Mathieu-Zhao subspace of $V$.

\begin{proof} Let $v\in V^0$, and let $\mathcal T=\{v(s_2)\cdots v(s_t)\textbf{1}|s_2,\dots,s_t\in\mathbb Z\}$.  Since $V^0$ is a subspace, we know
$$v(s_1)\cdots v(s_t)\textbf{1}\in V^0$$
when $t\geq 0$ for all $s_1,\dots,s_t\in\mathbb Z$.  Then for $0\neq v\in V^0$, we have $v\in\mathfrak{r}(V^i)$ if and only if there exists $m\geq 0$ such that
$$v(s_1)\cdots v(s_t)\textbf{1}=0$$
for all $t\geq m, s_1,\dots,s_t\in\mathbb Z$.  However, since $\text{Ann}_V(\mathcal T)$ is an ideal, and since $V$ is a simple VA, we have
$$\text{Ann}_V(\mathcal T)=\{0\},$$
thus $v=0$.  We repeat this argument for $j=1,\dots,k-1$ and note that for $v\in V^j\cap\mathfrak{r}(V^i)$, we must have $v=0$.

We now consider (by contradiction) the existence of $0\neq w:=\sum_{j=0}^{k-1}a_j\in \mathfrak{r}(V^i)$ such that $a_j\in V^j$.  Then by definition, we have $m\geq 0$ such that
$$w(s_1)\cdots w(s_t)\textbf{1}\in V^i$$
for all $t\geq m, s_1,\dots,s_t\in\mathbb Z$.  Denote
$$A:=w(s_1)\cdots w(s_m)\textbf{1}.$$
Since
$$w(p)A=\sum_{j=0}^{k-1}a_j(p)A\in V^i$$
for all $p\in\mathbb Z$, and since, for all $1\leq j\leq k-1$, we have $a_j(p)A\in V^{i+j}\neq V^i$, we can conclude that $a_j(p)A=0$ for all $p\in\mathbb Z$.  Then, again, since $\text{Ann}_V(\{A\})=\{0\}$, we conclude that $a_j=0$ for all $1\leq j\leq k-1$.  But this implies that $w=a_0\in V^0\cap\mathfrak{r}(V^i)=\{0\}$.  Thus $\mathfrak r(V^i)=\{0\}$ and $V^i$ is an MZ subspace of $V$.\end{proof}
\end{theorem}

\begin{theorem}[Main Result 2]
Consider $M_1(1)$.  Let $\{n_i|i\geq 1\}$ be a strictly increasing sequence such that $0\neq n_i\in\mathbb N$ for all $i\geq 1$ and $n_{i+1}-n_i\neq 1$ for infinitely many $i\geq 1$.  Let $M$ be the $\mathbb C$-subspace of $M_1(1)$ spanned by $\alpha(m_1)\cdots\alpha(m_{n_i})\textbf{1}$ over $\mathbb C$, where $m_1,\dots,m_{n_i}\in\mathbb N$.  Then the following are equivalent:
\begin{itemize}
\item[(1)] $\mathfrak{r}(M)=\{0\}$;
\item[(2)] $M$ is a Mathieu-Zhao subspace of $M_1(1)$;
\item[(3)] there exists no $d\in\mathbb{N}$ such that $md\in\{n_i|i\geq 1\}$ for all $m\geq 1$.
\end{itemize}

\begin{proof} $(1)\Rightarrow (2)$: See Corollary 36.
\\
$(2)\Rightarrow (3)$: Assume that there exists $d\geq 1$ such that $md\in\{n_i|i\geq 1\}$ for all $m\geq 1$.  Note that we must have $d>1$; otherwise, we would have $\{n_i\text{ }|\text{ }i\geq 1\}=\mathbb N$, a contradiction.  Then consider $d\geq 2$.  Then, by our assumption, we have
$$A:=\underbrace{\alpha(-1)\cdots\alpha(-1)}_{d\text{ times}}\textbf{1}\in M.$$
By Lemma 27, we know
$$A(q_1)\cdots A(q_m)\textbf{1}\in M,$$
for all $m\geq 1$ and $q_1,\dots,q_m\in\mathbb Z$, i.e., we know
$$A\in\mathfrak r(M).$$
Now, since $M$ is an MZ subspace of $M_1(1)$, for all $0\leq r\leq d-1$, there exists $N_r\geq 1$ such that
$$\big(A(q_1)\cdots A(q_m)\textbf{1}\big)(n)\underbrace{\alpha(-1)\cdots\alpha(-1)}_{r\text{ times}}\textbf{1}\in M$$
for all $m\geq N_r$ and $n\in\mathbb Z$.  This implies
$$\underbrace{\alpha(-1)\cdots\alpha(-1)}_{md+r\text{ times}}\textbf{1}\in M$$
for all $m\geq N_r$.  We consider $N:=\text{max}\{N_r|0\leq r\leq d-1\}$.  Then we have
$$\underbrace{\alpha(-1)\cdots\alpha(-1)}_{t\text{ times}}\textbf{1}\in M$$
for all $t\geq N$, a contradiction of our construction of $\{n_i|i\geq 1\}$.  We conclude that there must exist no such $d\in\mathbb N$.
\\
$(3)\Rightarrow (1)$: Assume that there exists $0\neq A\in\mathfrak r(M)$.  We may assume 
$$A\in\bigoplus_{i=0}^kM_1^i(1)$$
for some $k\geq 0$.  Then there exists $m\in\mathbb N$ such that
$$A(q_1)\cdots A(q_t)\textbf{1}\in M$$
for all $t\geq m$.  Let $c$ denote the maximum $\alpha$-length of $A\in M$.  Then
$$\underbrace{\alpha(-1)\cdots\alpha(-1)}_{ct\text{ times}}\textbf{1}\in M.$$
This means that $ct\in\{n_i|i\geq 1\}$ for all $t\geq m$, so let $d=cm$.  This implies that $md\in\{n_i|i\geq 1\}$ for all $m\geq 1$, a contradiction.  We conclude that $\mathfrak r(M)=\{0\}$.
\end{proof}
\end{theorem}

\begin{corollary}
Consider the eigenspace decomposition 
$$M_1(1)=\bigoplus_{l=0}^{k-1}M_1^{(l)}$$
under the automorphism $g$ of order $k$.  Then for any set $\{l_r|1\leq r\leq s\leq l_s\leq k-1\}$ of integers, the subspace
$$M:=\bigoplus_{r=1}^sM_1^{(l_r)}(1)$$
is a Mathieu-Zhao subspace of $M_1(1)$.

\begin{proof}First, note that since $0\notin\{l_r|1\leq r\leq s\leq k-1\}$, we may construct a sequence $\{n_i|i\geq 1\}$, with $n_1=l_1,\dots,n_s=l_s,n_{s+1}=l_1+k,\dots$, which satisfies that $n_{i+1}-n_i\neq 1$ for infinitely many $i\geq 1$.  Then we may apply Theorem 46 and consider condition $(3)$.  Suppose there exists $d\in\mathbb N$ such that $dm\in\{n_i|i\geq 1\}$ for all $m\geq 1$.  But this means that $dk\in\{n_i|i\geq 1\}$, a contradiction of our construction of $\{n_i|i\geq 1\}$.  We conclude that there exists no such $d\in\mathbb N$, thus $M$ is an MZ subspace of $M_1(1)$.\end{proof}
\end{corollary}

We have now succeeded in classifying an infinite set of nontrivial, non-ideal MZ subspaces of a (not necessarily simple) vertex algebra.  To conclude, we suggest next steps in our effort to survey MZ subspaces of vertex algebras.

\begin{center}
\textbf{Next Steps}
\end{center}

What Proposition 40 shows us is that we need not stray too far from associative algebras to learn something new.  After all, the MZ subspace structure is still in its infancy in terms of mathematical concepts, so pretty much anything we learn is new.  While there can be no wrong first step moving forward, it can be quite hard to tell what is too ambitious a goal at these early stages.  For example, studying the converse of Proposition 40 will prove quite difficult.  One would expect that such a statement would not hold, given how much more restrictive the MZ criteria turn out to be for VA's.  But with so much theory still to be developed and with literally infinitely many more examples (and possibly counterexamples) out there to find, it is still very hard to tell.

Of vital importance to the study of vertex (operator) algebras is (\cite{zhu1}).  Therein, Zhu introduced what has come to be known as the Zhu algebra of a vertex operator algebra.  This structure demonstrates that it is actually to be expected that we find such natural extensions from associative algebras to VA's.  The Zhu algebra of a given vertex algebra is defined as follows:

\begin{definition}
Let $V$ be a vertex operator algebra, and let $a\in V$ be homogeneous.  Define a bilinear operation
\[\begin{array}{l l l l}
*:&V\times V&\to& V,\\
&(a,b)&\mapsto&\text{Res}_z\Big(Y(a,z)\frac{(z+1)^{\text{deg }a}}{z}b\Big).
\end{array}\]
Denote the subspace of $V$ spanned by elements
$$\text{Res}_z\Big(Y(a,z)\frac{(z+1)^{\text{deg }a}}{z^2}b\Big)$$
by $O(V)$.  Then the quotient space $V/O(V):=A(V)$ is an associative algebra called the \textit{Zhu algebra of $V$}.
\end{definition}

\begin{theorem}
Let $V=M(1)$.  Then the Zhu algebra of $V$ is $$A(V)=V/O(V)\cong \mathbb C[x_i|i\geq 1].$$\qed
\end{theorem}

We continue to search for powerful qualifiers such as Zhao's so-called idempotent criterion for associative algebras (\cite{zhaoidempotent}) as follows:

\begin{theorem}
Let $\mathcal A$ be an algebraic algebra over a field $K$ and $M$ be a subspace of $\mathcal A$.  Then $M$ is a Mathieu-Zhao subspace of $\mathcal A$ if and only if for every idempotent $e\in M$, the ideal generated by $e$, denoted $(e)$, satisfies $(e)\subseteq M$.\qed
\end{theorem}

For a detailed proof of this powerful theorem, see (\cite{nieman, vanhove, zhaoidempotent}).  We use ideas in (\cite{localvoa}) to develop a notion of idempotents for vertex algebras, and we use this notion to attempt an analogue conjecture to Theorem 50 for vertex algebras.  The following definition, stated as a theorem in (\cite{localvoa}), gives a subspace on which hinges our definition of an idempotent of a vertex algebra.

\begin{definition}
Let
$$V=\bigoplus_{n\in\mathbb Z}V_n$$
be a $\mathbb Z$-graded vertex operator algebra.  Then we define the \textit{center of $V$} to be $$Z(V)=\{v\in V|Y(v,z)=v(-1)\}.$$
The set of idempotents $\{e_i\}\subseteq Z(V)$ is exactly the set of idempotents of $V$.
\end{definition}

A few more remarks from (\cite{localvoa}) allow us to more clearly state the following conjecture, which is analogous to the idempotent criterion for associative algebras:

\begin{conjecture}
Let $\{e_1,\dots,e_m\}$ be the unique set of primitive idempotents of a $\mathbb Z$-graded VOA,
$$V=\bigoplus_{n\in\mathbb Z}V_n.$$
We may decompose $V$ into ideals
$$V=V^1\oplus\cdots\oplus V^m$$
where
$$V^i=e_i(-1)V.$$
Then a subspace $M\subseteq V$ is a Mathieu-Zhao subspace of $V$ if and only if whenever $e_i\in M$, we have $V^i\subseteq M$.
\end{conjecture}

\begin{proposition}
Let $M\subseteq V$ be a Mathieu-Zhao subspace of a commutative vertex operator algebra $V$, and $\{e_1,\dots,e_m\}$ be the unique set of idempotents of $V$.  Then for each idempotent $e_i\in M$, we have $V^i\subseteq M$.

\begin{proof} Let $e_i\in M$.  Since $M$ is MZ, we understand that for all $w\in V$, there exists $k\geq 0$ such that for all $t\geq k$, and $n_1,\dots,n_t\in\mathbb Z$, we have
$$e_i(n_1)\cdots e_i(n_t)w\in M.$$
In particular, we have
$$\underbrace{e_i(-1)\cdots e_i(-1)}_{t\text{ times}}w=e_i(-1)w\in M\Leftrightarrow V^i\subseteq M.$$\end{proof}
\end{proposition}

We observe that $$Z(V)\cong \text{End}(V);$$
for a complete proof, see (\cite{localvoa}).  Should Conjecture 52 hold, we would significantly reduce the work of proving whether or not a given subspace of a VOA is MZ.  The ability to limit our search to a much more restricted set of operators could allow us to see, even more concretely, the connections between MZ subspaces and nilpotency.

\begin{conjecture}[LNED Conjecture for Vertex Algebras]
Let $K$ be a field of characteristic 0, $V$ be a vertex algebra over $K$, and $D$ be a locally nilpotent derivation of $V$.  Then for every ideal $I$ of $V$, the image $D(I)$ of $I$ under $D$ is a Mathieu-Zhao subspace of $V$.
\end{conjecture}

\newpage
\begin{center}
REFERENCES
\end{center}
\printbibliography[heading=none]

@article{zhu1,
author = {Yongchang Zhu},
title = {Modular Invariance of Characters of Vertex Operator Algebras},
journaltitle = {Journal of the American Mathematical Society},
year = {1996},
volume = {9},
number = {1},
pages = {237--302}
}

@thesis{vanhove,
author = {Loes C. van Hove},
title = {{Mathieu-Zhao} Subspaces},
type = {Master's Thesis},
institution = {Radboud University},
year = {2015}
}

@thesis{nieman,
author = {H. S. Nieman},
title = {Mathieu Subspaces of Univariate Polynomial Rings},
type = {Master's Thesis},
institution = {Radboud University},
year = {2012}
}

@misc{zhaoopen,
author = {Wenhua Zhao},
title = {Some Open Problems on Locally Finite or Locally Nilpotent Derivations and {$\varepsilon$}-Derivations},
year = {2017},
note = {arXiv:1701.05992v1}
}

@misc{zhaoimages,
author = {Wenhua Zhao},
title = {Images of Ideals under Derivations and {$\varepsilon$}-Derivations of Univariage Polynomial Algebras over a Field of Characteristic Zero},
year = {2017},
note = {arXiv:1701.06125v1}
}

@misc{dongintro,
author = {Chongying Dong},
title = {Introduction to Vertex Operator Algebras {I}},
year = {1994},
type = {Presented at Research Institute of Mathematical Sciences, Kyoto, Japan}
}

@article{localvoa,
author = {Chongying Dong and Geoffrey Mason},
title = {Local and Semilocal Vertex Operator Algebras},
year = {2004},
journaltitle = {Journal of Algebra},
volume = {280},
number = {1},
pages = {350--366}
}

@misc{borcherdstalk,
author = {Richard E. Borcherds},
title = {What is a Vertex Algebra?},
year = {1997},
type = {Presented at Mathematische Arbeitstagung, Max Planck Institute for Mathematics, Bonn, Germany}
}

@article{zhaogen,
author = {Wenhua Zhao},
title = {Generalizations of the Image Conjecture and the {Mathieu} Conjecture},
journal = {Journal of Pure and Applied Algebra},
year = {2010},
volume = {214},
number = {2010},
pages = {1200--1217}
}

@misc{zhaotalk,
author = {Wenhua Zhao},
title = {Recent Developments on {Mathieu} Subspaces},
year = {2017},
note = {Presented at Illinois State University Algebra Seminar}
}

@book{lepowskyli,
author = {James Lepowsky and Hiasheng Li},
title = {Introduction to Vertex Operator Algebras and their Representations},
year = {2004},
publisher = {Birkh\"auser}
}

@article{zhaoidempotent,
author = {Wenhua Zhao},
title = {Mathieu Subspaces of Associative Algebras},
journaltitle = {Journal of Algebra},
year = {2012},
volume = {350},
pages = {245--272}
}

@misc{essen,
author = {Arno van den Essen},
title = {The Amazing Image Conjecture},
year = {2010},
note = {arXiv:1006.5801v1}
}

@article{moonshine,
author = {Richard E. Borcherds},
title = {Monstrous Moonshine and Monstrous Lie Superalgebras},
journaltitle = {Inventiones Mathematicae},
year = {1992},
volume = {109},
number = {1},
pages = {405--444}
}

@article{pakovich,
author = {Fedor Pakovich},
title = {On Rational Functions Orthogonal to all Powers of a Given Rational Function on a Curve},
journaltitle = {Moscow Journal of Mathematics},
year = {2013},
volume = {13},
number = {4},
pages = {693--731}
}

@misc{intromathieu,
title = {Introduction to {Mathieu} Subspaces},
type = {Presented at International Short-School/Conference on Affine Algebraic Geometry and the Jacobian Conjecture, Chern Institute of Mathematics, Nankai University, Tianjin, China},
author = {Arno van den Essen},
year = {2014}
}

@article{fpyz,
author = {Jean-Pierre Francoise and Fedor Pakovich and Yousef Yomdin and Wenhua Zhao},
title = {Moment Vanishing Problem and Positivity: Some Examples},
journaltitle = {Bulletin des Sciences Math\'ematiques},
year = {2011},
volume = {13},
number = {1},
pages = {10--32}
}

@book{flm,
author = {Igor Frenkel and James Lepowski and Arne Meurman},
title = {Vertex Operator Algebras and the Monster},
year = {1988},
publisher = {Academic Press}
}

@article{borcherds1,
author = {Richard E. Borcherds},
title = {Vertex Algebras, {Kac-Moody} Algebras, and the Monster},
journaltitle = {Proceedings of the National Academy of Sciences, USA},
year = {1986},
pages = {3068--3071},
volume = {83}
}

\end{document}